\documentclass[10pt]{amsart}
\usepackage{amsfonts}
\usepackage{amssymb}
\usepackage{amsfonts}
\usepackage{tikz}

\newtheorem{thm1}{Theorem}[section]
\newtheorem{lem1}[thm1]{Lemma}
\newtheorem{rem1}[thm1]{Remark}

\newtheorem{cor1}[thm1]{Corollary}

\newtheorem{prop1}[thm1]{Proposition}
\newtheorem{ex1}[thm1]{Example}

\begin{document}

\title[Toric ideals and diagonal 2-minors]
{Toric ideals and diagonal 2-minors}
\author[A. Katsabekis]{Anargyros Katsabekis}
\address{Department of Mathematics, Mimar Sinan Fine Arts University, Istanbul, 34427, Turkey} \email{katsabek@aegean.gr}
\thanks{This work was carried out during the tenure of an ERCIM ``Alain Bensoussan" Fellowship Programme. The research leading to these results has received funding from the European Union Seventh Framework Programme (FP7/2007-2013) under grant agreement no 246016.}


\begin{abstract} Let $G$ be a simple graph on the vertex set $\{1,\ldots,n\}$ with $m$ edges. An algebraic object attached to $G$ is the ideal $P_{G}$ generated by diagonal 2-minors of an $n \times n$ matrix of variables. In this paper we prove that if $G$ is bipartite, then every initial ideal of $P_{G}$ is generated by squarefree monomials of degree at most $\left \lfloor{\frac{m+n+1}{2}} \right \rfloor$. Furthermore, we completely characterize all connected graphs $G$ for which $P_{G}$ is the toric ideal associated to a finite simple graph. Finally we compute in certain cases the universal Gr{\"o}bner basis of $P_{G}$.
\end{abstract}
\maketitle

\section{Introduction}

Let $X=(x_{ij})$ be an $n \times n$ matrix of variables and $R=K[x_{ij}|1 \leq i,j \leq n]$ be the polynomial ring in $n^2$ variables over a field $K$. For $1 \leq i<j \leq n$ we denote by $f_{ij}$ the diagonal 2-minor of $X$ given by the elements that stand at the intersection of the rows $i$, $j$ and the columns $i$, $j$, in other words, $f_{ij}=x_{ii}x_{jj}-x_{ij}x_{ji}$. Thus $f_{ij}$ is a binomial in $R$, namely a difference of two monomials.

Let $G$ be a simple graph on the vertex set $\{1,\ldots,n\}$ with the edge set $E(G)$. By a simple graph $G$ we mean an undirected graph without loops or multiple edges. Let $m$ be the cardinality of $E(G)$. We write $S$ for the polynomial ring $$S:=K[\{x_{ij},x_{ji}| \{i,j\} \in E(G)\} \cup \{x_{ii}\}_{1 \leq i \leq n}]$$ in $2m+n$ variables over $K$. Let $P_{G}$ be the ideal of $S$ generated by all the binomials $f_{ij}$, where $\{i,j\}$ is an edge of $G$.

Ideals generated by diagonal 2-minors were studied for the first time in \cite{EQ}. More specifically, the authors considered the ideal $P_{G}^{'}=P_{G} \cdot R$ and proved that it is a complete intersection prime ideal of height $m$. They also noticed that the set of all binomials $f_{ij}$, where $\{i,j\}$ is an edge of $G$, is the reduced Gr{\"o}bner basis of $P_{G}^{'}$ with respect to the reverse lexicographical order given by the natural ordering of variables $$x_{11}>x_{12}>\cdots>x_{1n}>x_{21}>x_{22}>\cdots>x_{2n}>
\cdots>x_{n1}>\cdots>x_{nn}.$$ By Theorem 1.2 in \cite{EQ}, the initial ideal of $P_{G}^{'}$ with respect to the lexicographic term order given by the natural order of indeterminates is generated by squarefree monomials of degree at most $4$.

The universal Gr{\"o}bner basis of an ideal $I$ generated by binomials, denoted by $U(I)$, is the union of all its reduced Gr{\"o}bner bases. Universal Gr{\"o}bner bases were used in integer programming and algebraic statistics to compute Markov bases of contingency tables and to study relations among conditional probabilities, see \cite{Mor} and \cite{Rap}. In general, it is not easy to describe the set $U(I)$ or even to estimate the degree of its elements. The case that $I$ is a toric ideal has been of particular interest, see for example \cite{BHP} and \cite{Tata-Tho}.

In this article we study the universal Gr{\"o}bner basis of $P_{G}$. Since the ideal $P_{G}$ is prime, it is also a toric ideal and therefore one can associate an integer program to it. The maximum degree of the elements in the universal Gr{\"o}bner basis of $P_{G}$ is an important measure for the complexity of the aforementioned integer program. It is worth to note that there exists an upper bound on the degrees of the elements of $U(P_{G})$ (see Theorem 4.7 in \cite{St}), however it is very far from being sharp.

In section 2 we show that if the graph $G$ is bipartite, then every initial ideal of $P_{G}$ is generated by squarefree monomials and also the universal Gr{\"o}bner basis consists only of the circuits of $P_{G}$. For a bipartite graph $G$, we additionally prove that the maximum degree of a binomial in the universal Gr{\"o}bner basis of $P_{G}$ is bounded above by $\left \lfloor{\frac{m+n+1}{2}} \right \rfloor$. Our bound is sharp, see Remark \ref{Upper bound sharp} and Remark \ref{Upper bound is sharp}.

An interesting problem is to determine when an ideal generated by binomials is the toric ideal $I_{H}$ associated to a finite graph $H$. In \cite{OH} H. Ohsugi and T. Hibi consider it for the class of ideals generated by adjacent 2-minors. In section 3 we study the above problem for ideals generated by diagonal 2-minors. More precisely, Theorem \ref{Th1 Toric} asserts that, given a connected graph $G$, there exists a finite simple graph $H$ such that $P_{G}=I_{H}$ if and only if $G$ has at most one cycle. Our approach produces new examples of non-bipartite graphs $H$ such that the toric ideal $I_{H}$ has a quadratic Gr{\"o}bner basis. Moreover, using Theorem \ref{Th1 Toric} and Theorem 3.4 in \cite{Tata-Tho}, we characterize in graph theoretical terms the elements of the universal Gr{\"o}bner basis of $P_{G}$, under the assumption that $G$ is a connected graph with at most one cycle.

\section{General results for ideals generated by diagonal $2$-minors}

In this section first we recall some basic facts about toric ideals associated to vector configurations. Next we associate to a simple graph $G$ on the vertex set $\{1,\ldots,n\}$ the vector configuration $A_{G}$, and let $N_{G}$ be the matrix with columns the vectors of $A_{G}$. It turns out (see Proposition \ref{Proposit}) that $P_{G}$ is the toric ideal $I_{A_{G}}$ associated to $A_{G}$. Also, we show that the matrix $N_{G}$ is totally unimodular if and only if the graph $G$ is bipartite. As a consequence, if $G$ is a bipartite graph, then every initial ideal of $P_{G}$ is generated by squarefree monomials. Moreover, for a bipartite graph $G$, the universal Gr{\"o}bner basis of $P_{G}$ is described by the circuits of $I_{A_{G}}$.

\subsection{Basics on toric ideals} \mbox{} \par

\medskip
Let $A=\{{\bf a}_{1},\ldots,{\bf a}_{s}\} \subset \mathbb{Z}^{r}$ be a vector configuration and let $K[y_{1},\ldots,y_{s}]$ be the polynomial ring in $s$ variables over a field $K$. As usual, we will denote by ${\bf y}^{\bf u}$ the monomial $y_{1}^{u_{1}} \cdots y_{s}^{u_{s}}$ of $K[y_{1},\ldots,y_{s}]$, with ${\bf u}=(u_1, \ldots, u_s) \in \mathbb{N}^s$. The {\em toric ideal} $I_{A}$ is the kernel of the $K$-algebra homomorphism $\phi: K[y_{1},\ldots,y_{s}] \rightarrow K[t_{1}^{\pm 1},\ldots,t_{r}^{\pm 1}]$ given by $$\phi(y_{i})= {\bf t}^{{\bf a}_i}=t_{1}^{a_{i,1}} \cdots t_{r}^{a_{i,r}}\ \ \textrm{for all} \ i=1,\ldots,s,$$ where ${\bf a}_{i}=(a_{i,1},\ldots,a_{i,r})$. Thus every toric ideal is prime. Furthermore, the ideal $I_{A}$ is generated by all the binomials ${\bf y}^{\bf u}-{\bf y}^{\bf v}$ such that $\phi({\bf y}^{\bf u})=\phi({\bf y}^{\bf v})$, see \cite{St}.

Let $N=(a_{i,j})$ be the $r \times s$ matrix with columns the vectors ${\bf a}_{1},\ldots,{\bf a}_{s}$. By Lemma 4.2 in \cite{St}, the height of $I_{A}$, denoted by ${\rm ht}(I_{A})$, is equal to $s-g$, where $g$ is the rank of the matrix $N$.

For a vector ${\bf u}=(u_{1},\ldots,u_{s}) \in \mathbb{Z}^{s}$, we let ${\rm supp}({\bf u})=\{i \in \{1,\ldots,s\}|u_{i} \neq 0\}$ be the support of ${\bf u}$. The support of a binomial $B={\bf y}^{{\bf u}}-{\bf y}^{{\bf v}}$ is ${\rm supp}(B)={\rm supp}({\bf u}) \cup {\rm supp}({\bf v})$. An irreducible binomial $B$ belonging to $I_{A}$ is called a {\em circuit} of $I_{A}$ if there exists no binomial $B' \in I_{A}$ such that ${\rm supp}(B') \subsetneqq  {\rm supp}(B)$. We shall denote by $\mathcal{C}(I_{A})$ the set of circuits of $I_{A}$. The set $\mathcal{C}(I_{A})$ can be computed easily, see \cite{Loera} or \cite{DES}. An irreducible binomial ${\bf y}^{{\bf u}}-{\bf y}^{{\bf v}}$ belonging to $I_{A}$ is called {\em primitive} if there is no other binomial ${\bf y}^{{\bf c}}-{\bf y}^{{\bf d}} \in I_{A}$ such that ${\bf y}^{{\bf c}}$ divides ${\bf y}^{{\bf u}}$ and ${\bf y}^{{\bf d}}$ divides ${\bf y}^{{\bf v}}$. The set of all primitive binomials is the Graver basis of $I_{A}$, denoted by ${\rm Gr}(I_{A})$. For any toric ideal $I_{A}$, we have that $\mathcal{C}(I_{A}) \subseteq U(I_{A}) \subseteq {\rm Gr}(I_{A})$, see \cite{St}.

Toric ideals associated to graphs serve as interesting examples of toric ideals. Let $H$ be a finite simple graph with vertices $V(H)=\{v_{1},\ldots,v_{r}\}$ and edges $E(H)=\{z_{1},\ldots,z_{s}\}$. The {\em incidence matrix} of $H$ is the $r \times s$ matrix $M_{H}=(b_{i,j})$ defined by $$b_{i,j}=\begin{cases} 1, & \textrm{if} \  v_{i} \ \textrm{is one of the vertices in} \ z_{j}\\
0, &  \textrm{otherwise}.
\end{cases}$$ Let $\mathcal{B}_{H}=\{{\bf b}_{1},\ldots,{\bf b}_{s}\}$ be the set of vectors in $\mathbb{Z}^{r}$, where ${\bf b}_{i}=(b_{i,1},\ldots,b_{i,r})$ for $1 \leq i \leq s$. With $I_{H}$ we denote the toric ideal $I_{\mathcal{B}_{H}}$ in $K[y_{1},\ldots,y_{s}]$. By Lemma 8.3.2 in \cite{Vila}, the rank of the matrix $M_{H}$ equals $r-c(H)$, where $c(H)$ is the number of connected components of $H$ which are bipartite. Thus the height of $I_{H}$ equals $s-r+c(H)$.

A {\em walk of length} $l$ of $H$ connecting $v_{i_{1}} \in V(H)$ with $v_{i_{l+1}} \in V(H)$ is a finite sequence of the form $$\Gamma=(\{v_{i_{1}},v_{i_{2}}\},\{v_{i_{2}},v_{i_{3}}\},\ldots,\{v_{i_{l}},v_{i_{l+1}}\})$$ with each $z_{i_j}=\{v_{i_j},v_{i_{j+1}}\} \in E(H)$, $1 \leq j \leq l$. In some cases we may also denote a walk only by vertices $(v_{i_1},v_{i_2},\ldots,v_{i_l},v_{i_{l+1}})$. An {\em even} (respectively {\em odd}) walk is a walk of even (respectively odd) length. The walk $\Gamma$ is {\em closed} if $v_{i_{1}}=v_{i_{l+1}}$. A {\em cycle} is a closed walk $(\{v_{i_{1}},v_{i_{2}}\},\{v_{i_{2}},v_{i_{3}}\},\ldots,\{v_{i_{l}},v_{i_{1}}\})$ with $l \geq 3$ and $v_{i_{j}} \neq v_{i_{k}}$, for every $1 \leq j<k \leq l$.

Given an even closed walk $\Gamma=(z_{i_1},\ldots,z_{i_{2l}})$ of $H$, we denote by $B_{\Gamma}$ the binomial $$B_{\Gamma}:=\prod_{j=1}^{l}y_{i_{2j-1}}-\prod_{j=1}^{l}y_{i_{2j}} \in I_{H}.$$ Moreover, Proposition 8.1.2 in \cite{Vila} guarantees that $I_{H}$ is generated by all the binomials $B_{\Gamma}$, where $\Gamma$ is an even closed walk of $H$.

\subsection{Ideals generated by diagonal 2-minors} \mbox{} \par
\medskip
Let $G$ be a simple graph on the vertex set $\{1,\ldots,n\}$ with the edge set $E(G)=\{z_{1},\ldots,z_{m}\}$. Let $\{{\bf e}_{1},\ldots,{\bf e}_{n},{\bf e}_{n+1},\ldots,{\bf e}_{n+m}\}$ be the canonical basis of $\mathbb{Z}^{n+m}$. Given an edge $z_{i}=\{i_{k},i_{l}\}$ of $G$, we consider the following two vectors: ${\bf a}_{\{i_{k}, i_{l}\}}:={\bf e}_{i_k}+{\bf e}_{i_l}-{\bf e}_{n+i}$ and ${\bf b}_{\{i_{k}, i_{l}\}}:={\bf e}_{n+i}$. Also, for every $1 \leq j \leq n$, we consider the vector ${\bf a}_{jj}:={\bf e}_{j}$. Let $A_{G}$ denote the set $$A_{G}:=\{{\bf a}_{\{i,j\}},{\bf b}_{\{i,j\}}| \{i,j\} \ \textrm{is an edge of} \ G\} \cup \{{\bf a}_{jj}|1 \leq j \leq n\}$$ consisting of $2m+n$ vectors. Let $N_{G}$ be the $(n+m) \times (2m+n)$ matrix whose columns are the vectors in $A_{G}$.
\begin{ex1} \label{Computing Example UGB} {\rm Let $G$ be the graph on the vertex set $\{1,\ldots,5\}$ with edges $z_{1}=\{1,2\}$, $z_{2}=\{2,3\}$, $z_{3}=\{3,4\}$, $z_{4}=\{1,4\}$, $z_{5}=\{1,5\}$ and $z_{6}=\{3,5\}$. The set $A_G$ consists of the following 17 vectors: ${\bf a}_{\{1,2\}}={\bf e}_{1}+{\bf e}_{2}-{\bf e}_{6}$, ${\bf b}_{\{1,2\}}={\bf e}_{6}$, ${\bf a}_{\{2,3\}}={\bf e}_{2}+{\bf e}_{3}-{\bf e}_{7}$, ${\bf b}_{\{2,3\}}={\bf e}_{7}$, ${\bf a}_{\{3,4\}}={\bf e}_{3}+{\bf e}_{4}-{\bf e}_{8}$, ${\bf b}_{\{3,4\}}={\bf e}_{8}$, ${\bf a}_{\{1,4\}}={\bf e}_{1}+{\bf e}_{4}-{\bf e}_{9},$ ${\bf b}_{\{1,4\}}={\bf e}_{9}$, ${\bf a}_{\{1,5\}}={\bf e}_{1}+{\bf e}_{5}-{\bf e}_{10}$, ${\bf b}_{\{1,5\}}={\bf e}_{10}$,  ${\bf a}_{\{3,5\}}={\bf e}_{3}+{\bf e}_{5}-{\bf e}_{11}$, ${\bf b}_{\{3,5\}}={\bf e}_{11}$, ${\bf a}_{11}={\bf e}_{1}$, ${\bf a}_{22}={\bf e}_{2}$, ${\bf a}_{33}={\bf e}_{3}$, ${\bf a}_{44}={\bf e}_{4}$, ${\bf a}_{55}={\bf e}_{5}$.\\ The toric ideal $I_{A_G}$ is the kernel of the $K$-algebra homomorphism $$\phi: S \rightarrow K[t_{1}^{\pm 1},\ldots,t_{11}^{\pm 1}]$$ given by $\phi(x_{12})=t_{1}t_{2}t_{6}^{-1}$, $\phi(x_{21})=t_{6}$, $\phi(x_{23})=t_{2}t_{3}t_{7}^{-1}$, $\phi(x_{32})=t_{7}$, $\phi(x_{34})=t_{3}t_{4}t_{8}^{-1}$, $\phi(x_{43})=t_{8}$,
$\phi(x_{14})=t_{1}t_{4}t_{9}^{-1}$, $\phi(x_{41})=t_{9}$,
$\phi(x_{15})=t_{1}t_{5}t_{10}^{-1}$, $\phi(x_{51})=t_{10}$,
$\phi(x_{35})=t_{3}t_{5}t_{11}^{-1}$, $\phi(x_{53})=t_{11}$, $\phi(x_{11})=t_{1}$, $\phi(x_{22})=t_{2}$, $\phi(x_{33})=t_{3}$, $\phi(x_{44})=t_{4}$, $\phi(x_{55})=t_{5}$. Actually $I_{A_G}$ is generated by the following $6$ binomials: $f_{12}$, $f_{23}$, $f_{34}$, $f_{14}$, $f_{15}$ and $f_{35}$. Thus $P_{G}=I_{A_G}$.}
\end{ex1}

Given a subgraph $W$ of $G$, we let $P_{W}$ be the ideal of the polynomial ring $$K[\{x_{ij},x_{ji}| \{i,j\} \in E(W)\} \cup \{x_{ii}\}_{1 \leq i \leq n}]$$ generated by all the binomials $f_{ij}$, where $\{i,j\}$ is an edge of $W$.

\begin{prop1} \label{Proposit1} (\cite{EQ}) The ideal $P_{W}$ is complete intersection of height $|E(W)|$.
\end{prop1}

\begin{prop1} \label{Proposit} Let $G$ be a simple graph on the vertex set $\{1,\ldots,n\}$ with $m$ edges, then the ideal $P_{G}$ coincides with the toric ideal $I_{A_{G}}$.
\end{prop1}
\noindent \textbf{Proof.} Given an edge $z_{i}=\{i_{k},i_{l}\}$ of $G$, we have that $${\bf a}_{\{i_{k}, i_{l}\}}+{\bf b}_{\{i_{k}, i_{l}\}}=({\bf e}_{i_k}+{\bf e}_{i_l}-{\bf e}_{n+i})+{\bf e}_{n+i}={\bf e}_{i_k}+{\bf e}_{i_l}={\bf a}_{i_{k}i_{k}}+{\bf a}_{i_{l}i_{l}},$$ so the binomial $f_{i_{k}i_{l}}$ belongs to $I_{A_G}$, and therefore $P_{G} \subseteq I_{A_{G}}$. Moreover the rank of $N_{G}$ is equal to $n+m$, thus ${\rm ht}(I_{A_{G}})=2m+n-(n+m)=m$. From Proposition \ref{Proposit1} we have that $P_{G}$ is a prime ideal of height $m$. Thus $P_{G}=I_{A_{G}}$. \hfill $\square$

\begin{rem1} \label{Rem Toric} {\rm There is a term order such that the initial ideal of $I_{A_{G}}$ is generated by squarefree quadratic monomials. By Corollary 8.9 in \cite{St}, the vector configuration $A_{G}$ has a unimodular regular triangulation.}
\end{rem1}

A matrix $M$ with ${\rm rank}(M)=d$ is called {\em unimodular} if all nonzero $d \times d$-minors of $M$ have the same absolute value. The matrix $M$ is called {\em totally unimodular} when every minor of $M$ is $0$ or $\pm 1$. It is well known (see for example \cite[Chapter 19]{Sch}) that the graph $G$ is bipartite if and only if its incidence matrix $M_{G}$ is totally unimodular.

\begin{thm1} \label{Unimodular} The matrix $N_{G}$ is totally unimodular if and only if the graph $G$ is bipartite.
\end{thm1}
\noindent \textbf{Proof.} ($\Rightarrow$) Suppose that $N_{G}$ is totally unimodular. Then also $M_{G}$ is totally unimodular since it is a submatrix of $N_{G}$, and therefore $G$ is bipartite.\\
($\Leftarrow$) Suppose that $G$ is bipartite. Since total unimodularity is preserved under the unit vectors, it is enough to consider the matrix $Q$ with columns the vectors ${\bf a}_{\{i_{k}, i_{l}\}}={\bf e}_{i_k}+{\bf e}_{i_l}-{\bf e}_{n+i}$, where $z_{i}=\{i_{k},i_{l}\}$ is an edge of $G$. It suffices to prove that $Q$ is totally unimodular. The incidence matrix $M_{G}$ is totally unimodular, since $G$ is bipartite, so also the transpose $M_{G}^{T}$ of $M_{G}$ is totally unimodular. Thus the $m \times (n+m)$ matrix $V=(M_{G}^{T}|R)$, where $R$ is the matrix with rows $-{\bf e}_{1},\ldots,-{\bf e}_{m}$, is totally unimodular, and therefore $Q=V^{T}$ is totally unimodular. \hfill $\square$

\begin{thm1} \label{UGB} Let $G$ be a simple bipartite graph on the vertex set $\{1,\ldots,n\}$, then the initial ideal of $P_{G}$ with respect to any term order is generated by squarefree monomials. Furthermore, the equality $U(P_{G})=\mathcal{C}(I_{A_{G}})$ holds.
\end{thm1}

\noindent \textbf{Proof.} We have, from Theorem \ref{Unimodular}, that the matrix $N_{G}$ is totally unimodular, and hence also unimodular. By Corollary 8.9 in \cite{St}, every initial ideal of $I_{A_{G}}$ is generated by squarefree monomials. Now, Proposition 8.11 in \cite{St} asserts that $\mathcal{C}(I_{A_{G}})={\rm Gr}(I_{A_{G}})$. Thus $U(I_{A_{G}})=\mathcal{C}(I_{A_{G}})$, and therefore $U(P_{G})=\mathcal{C}(I_{A_{G}})$ since $P_{G}=I_{A_{G}}$. \hfill $\square$

\begin{rem1} {\rm Let $G$ be a simple bipartite graph on the vertex set $\{1,\ldots,n\}$, then the matrix $N_{G}$ is unimodular. For every circuit $B={\bf x}^{\bf u}-{\bf x}^{\bf v} \in I_{A_{G}}$ we have, from the proof of Proposition 8.11 in \cite{St}, that both ${\bf x}^{\bf u}$ and ${\bf x}^{\bf v}$ are squarefree monomials.}
\end{rem1}

\begin{cor1} \label{IntegerPart} Let $G$ be a simple bipartite graph on the vertex set $\{1,\ldots,n\}$ with $m$ edges, then the maximum degree of a binomial in the universal Gr{\"o}bner basis of $P_{G}$ is less than or equal to $\left \lfloor{\frac{m+n+1}{2}} \right \rfloor$.
\end{cor1}

\noindent \textbf{Proof.} Given a circuit $B \in I_{A_{G}}$, we have, from Lemma 4.8 in \cite{St}, that ${\rm supp}(B)$ has at most $m+n+1$ elements. Since every binomial in $I_{A_{G}}$ is homogeneous with respect to the standard grading, we have that the degree of $B$ is less than or equal to $\left \lfloor{\frac{m+n+1}{2}} \right \rfloor$. By Theorem \ref{UGB}, the maximum degree of a binomial in the universal Gr{\"o}bner basis of $P_{G}$ is less than or equal to $\left \lfloor{\frac{m+n+1}{2}} \right \rfloor$. \hfill $\square$

\begin{ex1} {\rm We come back to Example \ref{Computing Example UGB}. The circuits of $I_{A_{G}}$ are the following: $$\mathcal{C}(I_{A_{G}})=\{f_{12}, f_{23}, f_{34}, f_{14}, f_{15}, f_{35}, x_{44}x_{35}x_{53}-x_{55}x_{34}x_{43},$$ $$x_{11}x_{35}x_{53}-x_{33}x_{15}x_{51}, x_{44}x_{15}x_{51}-x_{55}x_{14}x_{41}, x_{22}x_{15}x_{51}-x_{55}x_{12}x_{21},$$ $$x_{33}x_{14}x_{41}-x_{11}x_{34}x_{43}, x_{22}x_{14}x_{41}-x_{44}x_{12}x_{21}, x_{22}x_{34}x_{43}-x_{44}x_{23}x_{32},$$ $$x_{11}x_{23}x_{32}-x_{33}x_{12}x_{21}, x_{22}x_{35}x_{53}-x_{55}x_{23}x_{32}, x_{14}x_{41}x_{35}x_{53}-x_{34}x_{43}x_{15}x_{51},$$ $$x_{14}x_{41}x_{35}x_{53}-x_{33}x_{44}x_{15}x_{51}, x_{14}x_{41}x_{35}x_{53}-x_{11}x_{55}x_{34}x_{43},$$ $$x_{11}x_{44}x_{35}x_{53}-x_{34}x_{43}x_{15}x_{51}, x_{12}x_{21}x_{35}x_{53}-x_{23}x_{32}x_{15}x_{51},$$ $$x_{11}x_{22}x_{35}x_{53}-x_{23}x_{32}x_{15}x_{51}, x_{12}x_{21}x_{35}x_{53}-x_{22}x_{33}x_{15}x_{51},$$ $$x_{12}x_{21}x_{35}x_{53}-x_{11}x_{55}x_{23}x_{32}, x_{34}x_{43}x_{15}x_{51}-x_{33}x_{55}x_{14}x_{41},$$ $$x_{23}x_{32}x_{15}x_{51}-x_{33}x_{55}x_{12}x_{21}, x_{23}x_{32}x_{14}x_{41}-x_{12}x_{21}x_{34}x_{43},$$ $$x_{23}x_{32}x_{14}x_{41}-x_{11}x_{22}x_{34}x_{43}, x_{22}x_{33}x_{14}x_{41}-x_{12}x_{21}x_{34}x_{43},$$ $$x_{23}x_{32}x_{14}x_{41}-x_{33}x_{44}x_{12}x_{21}, x_{12}x_{21}x_{34}x_{43}-x_{11}x_{44}x_{23}x_{32},$$ $$x_{22}x_{14}x_{41}x_{35}x_{53}-x_{44}x_{23}x_{32}x_{15}x_{51}, x_{22}x_{14}x_{41}x_{35}x_{53}-x_{55}x_{12}x_{21}x_{34}x_{43},$$ $$ x_{44}x_{12}x_{21}x_{35}x_{53}-x_{55}x_{23}x_{32}x_{14}x_{41}, x_{44}x_{12}x_{21}x_{35}x_{53}-x_{22}x_{34}x_{43}x_{15}x_{51},$$ $$x_{22}x_{34}x_{43}x_{15}x_{51}-x_{55}x_{23}x_{32}x_{14}x_{41}, x_{44}x_{23}x_{32}x_{15}x_{51}-x_{55}x_{12}x_{21}x_{34}x_{43}\}.$$ Notice that $G$ is a bipartite graph which has more than one cycles. By Theorem \ref{UGB}, the universal Gr{\"o}bner basis of $P_{G}$ consists of the above 36 binomials.}
\end{ex1}

\section{Classification of all graphs $G$ such that the equality $P_{G}=I_{H}$ holds}

Let $G$ be a simple graph on the vertex set $\{1,\ldots,n\}$. The next proposition gives a necessary condition for the equality $P_{G}=I_{H}$, where $H$ is a finite simple graph.

\begin{prop1} \label{Basic Proposition} Let $G$ be a connected simple graph on the vertex set $\{1,\ldots,n\}$ with $m$ edges. If there exists a finite simple graph $H$ such that $P_{G}=I_{H}$, then $G$ has at most one cycle.
\end{prop1}

\noindent \textbf{Proof.} Without loss of generality we can assume that $H$ has exactly $2m+n$ edges and also that $H$ contains no isolated vertices. Recall that $P_{G}=I_{A_{G}}$, where $A_{G}=\{{\bf a}_{\{i,j\}},{\bf b}_{\{i,j\}}| \{i,j\} \in E(G)\} \cup \{{\bf a}_{jj}|1 \leq j \leq n\}$. Let $$\mathcal{B}_{H}=\{{\bf b}_{ij},{\bf b}_{ji}| \{i,j\} \in E(G)\} \cup \{{\bf b}_{11},\ldots,{\bf b}_{nn}\}$$ be the set of columns of the incidence matrix $M_{H}$ of $H$. Every column of $M_{H}$ has exactly two nonzero entries and both of them are equal to $1$. In particular, every ${\bf b}_{ii}$, $1 \leq i \leq n$, has exactly two nonzero entries, so the cardinality of the set $\{{\rm supp}({\bf b}_{11}),\ldots,{\rm supp}({\bf b}_{nn})\}$ is at most $2n$. Given an edge $\{i,j\}$ of $G$, we have that $f_{ij}=x_{ii}x_{jj}-x_{ij}x_{ji} \in P_{G}$, so $f_{ij} \in I_{H}$ and therefore ${\bf b}_{ij}={\bf b}_{ii}+{\bf b}_{jj}-{\bf b}_{ji}$. If ${\rm supp}({\bf b}_{ji}) \cap {\rm supp}({\bf b}_{ii})=\emptyset$ or ${\rm supp}({\bf b}_{ji}) \cap {\rm supp}({\bf b}_{jj})=\emptyset$, then the vector ${\bf b}_{ij}$ has at least one entry which is equal to $-1$, a contradiction. Thus the cardinality of the set $$\mathcal{R}=\{{\rm supp}({\bf b}_{ij}),{\rm supp}({\bf b}_{ji})| \{i,j\} \in E(G)\} \cup \{{\rm supp}({\bf b}_{11}),\ldots,{\rm supp}({\bf b}_{nn})\}$$ is at most $2n$. Let $d$ be the number of vertices of $H$, then $d$ is less than or equal to $2n$. Let $c(H)$ be the number of connected components of $H$ which are bipartite. By Proposition \ref{Proposit1} we have that ${\rm ht}(P_{G})=m$. Using the equality ${\rm ht}(P_{G})={\rm ht}(I_{H})$, we deduce that $m=2m+n-d+c(H)$, and therefore $m+n+c(H)=d \leq 2n$. So $m \leq m+c(H) \leq n$, while $n \leq m+1$ since $G$ is connected. Thus $n \in \{m,m+1\}$, and therefore $G$ has at most one cycle. Note that $G$ is a tree when $n=m+1$, while $G$ has exactly one cycle in the case that $n=m$. \hfill $\square$\\

Let $W$ be a connected simple graph with $k$ vertices and $l$ edges. Consider two simple graphs $W_{1}$, $W_{2}$ with $V(W_{1})=\{p_{1},\ldots,p_{k}\}$ and $V(W_{2})=\{q_{1},\ldots,q_{k}\}$ such that $V(W_{1}) \cap V(W_{2})=\emptyset$. Suppose that both $W_{1}$ and $W_{2}$ are isomorphic to $W$. More precisely, we assume that $\tau: V(W) \longrightarrow V(W_{1})$, $i \mapsto p_{i}$, is a bijection such that $$\{i,j\} \in E(W) \Longleftrightarrow \{p_{i},p_{j}\} \in E(W_{1}),$$ and also $\psi: V(W) \longrightarrow V(W_{2}), i \mapsto q_{i}$, is a bijection such that $$\{i,j\} \in E(W) \Longleftrightarrow \{q_{i},q_{j}\} \in E(W_{2}).$$ We will define a new graph $W^{\star}$ with $2k$ vertices and $k+2l$ edges as follows. The vertex set of $W^{\star}$ is $V(W_{1}) \cup V(W_{2})$. Also both $W_{1}$ and $W_{2}$ are subgraphs of $W^{\star}$, so each one of the edges of $W_{1}$, $W_{2}$ is also an edge of $W^{\star}$. Finally, for every vertex $i$ of $W$, we let $\{p_{i},q_{i}\}$ be an edge of $W^{\star}$. Note that $W^{\star}$ is a connected graph. The graph $W^{\star}$ is commonly known as a prism of $W$.

\begin{rem1} {\rm Let $G$ be a simple graph on the vertex set $\{1,\ldots,n\}$ and let $W$ be a connected subgraph of $G$. Given an edge $\{i,j\}$ of $W$, we have that $$\Gamma=(d_{ii}=\{p_{i},q_{i}\},d_{ji}=\{q_{i},q_{j}\},d_{jj}=\{p_{j},q_{j}\},
d_{ij}=\{p_{i},p_{j}\})$$ is an even cycle of $W^{\star}$, and therefore $x_{ii}x_{jj}-x_{ij}x_{ji}=B_{\Gamma} \in I_{W^{\star}}$. Thus $P_{W} \subseteq I_{W^{\star}}$.}
\end{rem1}

\begin{lem1} \label{Tree} Let $G$ be a simple graph on the vertex set $\{1,\ldots,n\}$ and $W$ be a connected subgraph of $G$. If $W$ is either a tree or a non-bipartite graph with exactly one cycle, then $P_{W}=I_{W^{\star}}$ where $W^{\star}$ is the graph constructed above.
\end{lem1}

\noindent \textbf{Proof.} Let $W$ be a tree with $k$ vertices and $l=k-1$ edges. By Proposition \ref{Proposit1} we have that ${\rm ht}(P_{W})=l$. If the graph $W^{\star}$ is not bipartite, then $${\rm ht}(I_{W^{\star}})=(k+2l)-2k=2l-k=2l-(l+1)=l-1<l,$$ and therefore ${\rm ht}(I_{W^{\star}})<{\rm ht}(P_{W})$ a contradiction to the fact that $P_{W} \subseteq I_{W^{\star}}$. Thus $W^{\star}$ is bipartite, so ${\rm ht}(I_{W^{\star}})=(k+2l)-2k+1=l$, and therefore $P_{W}=I_{W^{\star}}$.

Let $W$ be a non-bipartite graph with exactly one cycle. Let $r$ and $s$ be the number of vertices and edges, respectively, of $W$, then $r=s$. We have that the graph $W^{\star}$ is not bipartite, and therefore $${\rm ht}(I_{W^{\star}})=(r+2s)-2r=2s-r=s={\rm ht}(P_{W}).$$ Since $P_{W} \subseteq I_{W^{\star}}$, we obtain the equality $P_{W}=I_{W^{\star}}$. \hfill $\square$

\begin{rem1} {\rm Given an even cycle $C=(\{1,2\},\{2,3\},\ldots,\{k,1\})$ of a simple graph $G$, we have that the graph $C^{\star}$ is bipartite, since $$\{p_{1},p_{3},p_{5},\ldots,p_{k-1},q_{2},q_{4},\ldots,q_{k}\} \cup \{p_{2},p_{4},\ldots,p_{k},q_{1},q_{3},q_{5},\ldots,q_{k-1}\}$$ is a partition of its vertices. Thus ${\rm ht}(I_{C^{\star}})=3k-2k+1=k+1 \neq k={\rm ht}(P_{C})$, and therefore $P_{C} \neq I_{C^{\star}}$.}
\end{rem1}

Let $C=(\{i_{1},i_{2}\},\{i_{2},i_{3}\}, \ldots,\{i_{k},i_{1}\})$ be an even cycle of a graph $G$ of length $k \geq 4$. Consider two simple graphs $C_{1}$, $C_{2}$ with $V(C_{1})=\{v_{1},\ldots,v_{k}\}$ and $V(C_{2})=\{w_{1},\ldots,w_{k}\}$ such that $V(C_{1}) \cap V(C_{2})=\emptyset$. Suppose that both $C_{1}$ and $C_{2}$ are isomorphic to the path $Y=(\{i_{1},i_{2}\},\{i_{2},i_{3}\}, \ldots,$ $\{i_{k-1},i_{k}\})$. More precisely, we assume that $\tau: V(Y) \longrightarrow V(C_{1})$, $r \mapsto v_{r}$, is a bijection such that $$\{r,s\} \in E(Y) \Longleftrightarrow z_{rs}:=\{v_{r},v_{s}\} \in E(C_{1}),$$ and also $\psi: V(Y) \longrightarrow V(C_{2}), r \mapsto w_{r}$, is a bijection such that $$\{r,s\} \in E(Y) \Longleftrightarrow z_{sr}:=\{w_{r},w_{s}\} \in E(C_{2}).$$ We will define a new graph $\overline{C}$ with $2k$ vertices and $3k$ edges as follows. The vertex set of $\overline{C}$ is $V(C_{1}) \cup V(C_{2})$. Also both $C_{1}$ and $C_{2}$ are subgraphs of $\overline{C}$. Additionally, for $r=i_{1}$ and $s=i_{k}$, we let $z_{rs}:=\{v_{r},w_{s}\}$ and $z_{sr}:=\{v_{s},w_{r}\}$ be edges of $\overline{C}$. Finally, for every vertex $r$ of $C$, we let $z_{rr}:=\{v_{r},w_{r}\}$ be an edge of $\overline{C}$. Note that $\overline{C}$ is a connected graph. The graph $\overline{C}$ is commonly known as a twisted prism of $C$. The next lemma asserts that $P_{C}=I_{\overline{C}}$.

\begin{lem1} \label{CycleToric} Let $G$ be a simple graph on the vertex set $\{1,\ldots,n\}$ and $C$ be an even cycle of $G$ of length $k \geq 4$. Then $P_{C}=I_{\overline{C}}$ for the graph $\overline{C}$ constructed above.
\end{lem1}

\noindent \textbf{Proof.} Let $C=(\{i_{1},i_{2}\},\{i_{2},i_{3}\}, \ldots,\{i_{k},i_{1}\})$, then the graph $\overline{C}$ is not bipartite since $\Gamma=(v_{i_1},w_{i_k},v_{i_k},v_{i_{k-1}},\ldots,v_{i_2},v_{i_1})$ is an odd cycle of length $k+1$. Thus ${\rm ht}(I_{\overline{C}})=3k-2k=k$. For every edge $\{r,s\}$ of $Y$ we have that $\Gamma=(z_{rr},z_{sr},z_{ss},z_{rs})$ is an even cycle of $\overline{C}$, so $x_{rr}x_{ss}-x_{rs}x_{sr}=B_{\Gamma} \in I_{\overline{C}}$. Moreover, for $r=i_{1}$ and $s=i_{k}$, we have that also $\Gamma=(z_{rr},z_{sr},z_{ss},z_{rs})$ is an even cycle of $\overline{C}$, so $x_{rr}x_{ss}-x_{rs}x_{sr} \in I_{\overline{C}}$, and therefore $P_{C} \subseteq I_{\overline{C}}$. Thus the equality $P_{C}=I_{\overline{C}}$ holds, since, from Proposition \ref{Proposit1}, ${\rm ht}(P_{C})=k$. \hfill $\square$

\begin{rem1} \label{Odd Cycles Basic} {\rm Given an even cycle $C=(\{1,2\},\{2,3\},\ldots,\{k,1\})$ of length $k \geq 4$ of $G$, we have that the graph $\overline{C}$ contains at least one even cycle of length $2k$, namely $$(v_{1}, v_{2},\ldots, v_{k-1}, v_{k}, w_{k}, w_{k-1},\ldots,w_{1}, v_{1}).$$}
\end{rem1}

Let $G_{1}=(V(G_{1}),E(G_{1}))$, $G_{2}=(V(G_{2}),E(G_{2}))$ be graphs such that $G_{1} \cap G_{2}$ is a complete graph. The new graph $G=G_{1} \bigoplus G_{2}$ with the vertex set $V(G)=V(G_{1}) \cup V(G_{2})$ and edge set $E(G)=E(G_{1}) \cup E(G_{2})$ is called the {\em clique sum} of $G_{1}$ and $G_{2}$ in $G_{1} \cap G_{2}$. If the cardinality of $V(G_{1}) \cap V(G_{2})$ is $k+1$, then this operation is called a $k$-clique sum of the graphs $G_{1}$ and $G_{2}$. We write $G=G_{1} \bigoplus_{\widehat{v}} G_{2}$ to indicate that $G$ is the clique sum of $G_{1}$ and $G_{2}$ and that $V(G_{1}) \cap V(G_{2})=\widehat{v}$.

\begin{ex1} {\rm Let $G$ be the graph on the vertex set $\{1,\ldots,6\}$ consisting of the even cycle $C=(\{1,2\},\{2,3\},\{3,4\},\{1,4\})$ and the tree $T$ with edges $\{1,5\}$ and $\{1,6\}$. Let $\overline{C}$ be the graph on the vertex set $\{v_{1}=7,v_{2}=8,v_{3}=9,v_{4}=10\} \cup \{w_{1}=12,w_{2}=13,w_{3}=15,w_{4}=17\}$ consisting of the edges $z_{12}=\{7,8\}$, $z_{23}=\{8,9\}$, $z_{34}=\{9,10\}$, $z_{14}=\{7,17\}$, $z_{21}=\{12,13\}$, $z_{32}=\{13,15\}$, $z_{43}=\{15,17\}$, $z_{41}=\{10,12\}$, $z_{11}=\{7,12\}$, $z_{22}=\{8,13\}$, $z_{33}=\{9,15\}$, $z_{44}=\{10,17\}$.\\ Also consider the graph $T^{\star}$ on the vertex set $\{v_{1}=7,p_{5}=18,p_{6}=19\} \cup \{w_{1}=12,q_{5}=20,q_{6}=21\}$ consisting of the edges $z_{15}=\{7,18\}$, $z_{16}=\{7,19\}$, $z_{51}=\{12,20\}$, $z_{61}=\{12,21\}$, $z_{11}$, $z_{55}=\{18,20\}$, $z_{66}=\{19,21\}$. Notice that $\overline{C} \cap T^{\star}$ is the graph on the vertex set $\widehat{v}=\{7,12\}$ consisting of the edge $z_{11}$. The 1-clique sum $H$ of the graphs $\overline{C}$ and $T^{\star}$ is drawn in Figure 1.
\begin{center}
\begin{tikzpicture}
  [scale=.8,auto=left,every node/.style={circle,fill=blue!20}]
   \node (n14) at (1,1) {$17$};
   \node (n18) at (9,1) {$20$};
   \node (n17) at (9,7) {$18$};
   \node (n16) at (11,3) {$21$};
   \node (n15) at (11,5) {$19$};
   \node (n13) at (3,3) {$15$};
   \node (n11) at (7,3) {$12$};
   \node (n9) at (3,5) {$9$};
   \node (n12) at (5,3) {$13$};
    \node (n10) at (1,7) {$10$};
   \node (n8) at (5,5) {$8$};
   \node (n7) at (7,5) {$7$};

  \foreach \from/\to in {n7/n8,n7/n11,n7/n15,n7/n17,n8/n9,n8/n12,n9/n10,n9/n13,
  n10/n14,n13/n14,n10/n11,n7/n14,n11/n12,n12/n13,n11/n16,
  n11/n18,n15/n16,n17/n18}
    \draw (\from) -- (\to);

\end{tikzpicture}
\end{center}
\begin{center}
Figure 1.
\end{center}
It is easy to see that $P_{G} \subseteq I_{H}$. Moreover, ${\rm ht}(P_{G})=6$, and also ${\rm ht}(I_{H})=18-12=6$ since the graph $H$ is not bipartite. Thus $P_{G}=I_{H}$, so $\{f_{12},f_{23},f_{34},f_{14},f_{15},f_{16}\}$ is a quadratic Gr\"{o}bner basis for $I_{H}$ with respect to the reverse lexicographic term order given by $$x_{11}>x_{12}>x_{14}>x_{15}>x_{16}>x_{21}>x_{22}>x_{23}> x_{32}>x_{33}>
x_{34}>$$ $$x_{41}>x_{43}>x_{44}>x_{51}>x_{55}>x_{61}>x_{66}.$$}
\end{ex1}

Let $W$ be a subgraph of $G$. Given an ideal $J$ in the polynomial ring $$K[\{x_{ij},x_{ji}| \{i,j\} \in E(W)\} \cup \{x_{ii}\}_{1 \leq i \leq n}],$$ we write $J^{e}$ for the ideal $J^{e}:=J \cdot S$.

The following theorem determines all connected graphs $G$ such that the ideal $P_{G}$ is of the form $I_{H}$, for a finite simple graph $H$.

\begin{thm1} \label{Th1 Toric}  Let $G$ be a connected simple graph on the vertex set $\{1,\ldots,n\}$. Then, there exists a finite simple graph $H$ such that $P_{G}=I_{H}$ if and only if $G$ has at most one cycle.
\end{thm1}

\noindent \textbf{Proof.} If there exists a finite simple graph $H$ such that $P_{G}=I_{H}$, then we have, from Proposition \ref{Basic Proposition}, that $G$ has at most one cycle. Conversely if either $G$ has no cycle or $G$ is non-bipartite with exactly one cycle, then Lemma \ref{Tree} asserts that $P_{G}=I_{G^{\star}}$. Thus it is enough to consider the case that $G$ is bipartite with exactly one cycle. Let $C$ be the unique cycle of $G$ and suppose that it has length $k \geq 4$, where $k$ is even. The graph $G$ can be written as the $0$-clique sum of the cycle $C$ and some trees. More precisely, we have that $$G=C \bigoplus_{i_{1}} T_{1} \bigoplus_{i_{2}}T_{2} \bigoplus_{i_{3}} \cdots \bigoplus_{i_{s}}T_{s},$$ for some vertices $i_{1},\ldots,i_{s}$ of $C$. By Lemma \ref{CycleToric}, there is a connected graph $\overline{C}$ such that $P_{C}=I_{\overline{C}}$. Moreover, $\overline{C}$ has exactly $2k$ vertices and $3k$ edges. We shall denote by $\widehat{v_j}$ the set of vertices of the edge $z_{i_{j}i_{j}}$, $1 \leq j \leq s$. By Lemma \ref{Tree}, for each $1 \leq j \leq s$ there exists a connected graph $T_{j}^{\star}$ such that $P_{T_{j}}=I_{T_{j}^{\star}}$. Without loss of generality we can assume that $V(T_{j}^{\star}) \cap V(\overline{C})=\widehat{v_j}$, for every $1 \leq j \leq s$, and also $V(T_{j}^{\star}) \cap V(T_{k}^{\star})=\emptyset$, for every $j \neq k$. Let us suppose that the tree $T_{j}$, $1 \leq j \leq s$, has $g_{j}$ edges, then $T_{j}^{\star}$ has $2g_{j}+2$ vertices and $3g_{j}+1$ edges. Let
$$H=\overline{C} \bigoplus_{\widehat{v_1}}T_{1}^{\star} \bigoplus_{\widehat{v_2}}T_{2}^{\star} \bigoplus_{\widehat{v_3}} \cdots \bigoplus_{\widehat{v_s}}T_{s}^{\star}$$ be the 1-clique sum of the graphs $\overline{C}$ and $T_{j}^{\star}$, $1 \leq j \leq s$. We have that $P_{G}=(P_{C})^{e}+(P_{T_{1}})^{e}+\cdots+(P_{T_{s}})^{e}$, so  $P_{G}=(I_{\overline{C}})^{e}+(I_{T_{1}^{\star}})^{e}+\cdots+
(I_{T_{s}^{\star}})^{e}$, and therefore $P_{G} \subseteq I_{H}$. Notice that, from Proposition \ref{Proposit1}, ${\rm ht}(P_{G})=k+g_{1}+\cdots+g_{s}$, and also $${\rm ht}(I_{H})=(3k+3g_{1}+\cdots+3g_{s})-(2k+2g_{1}+\cdots+2g_{s})=k+
g_{1}+\cdots+g_{s}.$$ Consequently $P_{G}=I_{H}$. \hfill $\square$

\begin{rem1} {\rm Let $G$ be a connected simple graph on the vertex set $\{1,\ldots,n\}$ with exactly one cycle. Let $W$ be either the graph $H$ constructed in the proof of Theorem \ref{Th1 Toric} (when $G$ is bipartite) or the graph $G^{\star}$ when $G$ is non-bipartite. In both cases the graph $W$ is non-bipartite. The toric ideal $I_{W}$ has a quadratic Gr\"{o}bner basis with respect to the reverse lexicographic term order on $S$ induced by the natural ordering of variables $$x_{11}>x_{12}>\cdots>x_{1n}>x_{21}>x_{22}>\cdots>x_{2n}>\cdots>x_{n1}>\cdots>x_{nn}.$$}
\end{rem1}

\begin{ex1} \label{Graver} {\rm Consider the graph $G$ on the vertex set $\{1,\ldots,4\}$ with edges $d_{12}=\{1,2\}$, $d_{23}=\{2,3\}$, $d_{13}=\{1,3\}$ and $d_{14}=\{1,4\}$. The graph $G$ has exactly one odd cycle, namely $\Gamma=(d_{12}, d_{23}, d_{13})$. We let $G^{\star}$ be the graph with vertices $\{1,2,3,4\} \cup \{5,6,7,8\}$ and 12 edges, namely the 4 edges of $G$ and the edges $d_{21}=\{5,6\}$, $d_{32}=\{6,7\}$, $d_{31}=\{5,7\}$, $d_{41}=\{5,8\}$, $d_{11}=\{1,5\}$, $d_{22}=\{2,6\}$, $d_{33}=\{3,7\}$, $d_{44}=\{4,8\}$. By Lemma \ref{Tree}, the equality $P_{G}=I_{G^{\star}}$ holds. Thus the set $\{f_{12},f_{23},f_{13},f_{14}\}$ constitutes a quadratic Gr\"{o}bner basis for the toric ideal $I_{G^{\star}}$ with respect to the reverse lexicographic term order given by $$x_{11}>x_{12}>x_{13}>x_{14}>x_{21}>x_{22}>x_{23}>
x_{31}>x_{32}>x_{33}>x_{41}>x_{44}.$$ Using Algorithm 7.2 in \cite{St}, we determine the Graver basis of $I_{G^{\star}}$ which consists of the following 16 binomials: $$f_{12}, f_{23},f_{13}, f_{14}, x_{33}x_{14}x_{41}-x_{44}x_{13}x_{31}, x_{22}x_{14}x_{41}-x_{44}x_{12}x_{21},$$ $$x_{22}x_{13}x_{31}-x_{11}x_{23}x_{32}, x_{22}x_{13}x_{31}-x_{33}x_{12}x_{21}, x_{11}x_{23}x_{32}-x_{33}x_{12}x_{21},$$ $$x_{23}x_{32}x_{14}x_{41}-x_{22}x_{44}x_{13}x_{31}, x_{23}x_{32}x_{14}x_{41}-x_{33}x_{44}x_{12}x_{21}, x_{33}^{2}x_{12}x_{21}-x_{23}x_{32}x_{13}x_{31},$$ $$x_{11}^{2}x_{23}x_{32}-x_{12}x_{21}x_{13}x_{31}, x_{22}^{2}x_{13}x_{31}-x_{12}x_{21}x_{23}x_{32},$$ $$x_{11}x_{23}x_{32}x_{14}x_{41}-x_{44}x_{12}x_{21}x_{13}x_{31}, x_{12}x_{21}x_{13}x_{31}x_{44}^{2}-x_{23}x_{32}x_{14}^{2}x_{41}^{2}.$$ Notice that ${\rm Gr}(I_{G^{\star}}) \neq \mathcal{C}(I_{G^{\star}})$, since the binomial $B=x_{11}x_{23}x_{32}x_{14}x_{41}-x_{44}x_{12}x_{21}x_{13}x_{31}$ is primitive and not a circuit of $I_{G^{\star}}$. Thus $U(P_{G}) \neq \mathcal{C}(I_{G^{\star}})$.}
\end{ex1}

Let $G$ be a connected, simple and bipartite graph with exactly one cycle $C$. Consider the 1-clique sum $H$  of the graphs $\overline{C}$ and $T_{j}^{\star}$, $1 \leq j \leq s$, which appeared in the proof of Theorem \ref{Th1 Toric}. Note that $I_{A_{G}}=I_{H}$, since $P_{G}=I_{H}$. By Theorem \ref{UGB}, the equality $U(P_{G})=\mathcal{C}(I_{H})$ holds. Theorem \ref{UGB Exactlyone} describes all elements of $U(P_{G})$. In order to prove this theorem, we will use the following result:

\begin{thm1} \label{Circuits of a graph} (\cite{Vila}) Let $W$ be a finite, connected and simple graph. Then a binomial $B$ is a circuit of $I_{W}$ if and only if $B=B_{\Gamma}$ where $\Gamma$ is one of the following even closed walks of $W$: \begin{enumerate} \item $\Gamma$ is an even cycle.  \item $\Gamma=(C_{1},C_{2})$, where $C_{1}$ and $C_{2}$ are odd cycles of $W$ having exactly one common vertex.  \item $\Gamma=(C_{1},e_{1},\ldots,e_{r},C_{2},e_{r},\ldots,e_{1})$, where $C_{1}$ and $C_{2}$ are odd cycles of $W$ having no common vertex and where $(e_{1},\ldots,e_{r})$ is a path of $W$ which connects a vertex of $C_{1}$ and a vertex of $C_{2}$.
\end{enumerate}
\end{thm1}

\begin{thm1} \label{UGB Exactlyone} Let $G$ be a connected simple graph on the vertex set $\{1,\ldots,n\}$. Suppose that $G$ is bipartite with exactly one cycle. Then the universal Gr\"{o}bner basis of $P_{G}$ is given by the set $$U(P_{G})=\{B_{\Gamma}| \Gamma \ \textrm{is an even cycle of} \ H\}.$$
\end{thm1}

\noindent \textbf{Proof.} Let $C=(\{1,2\},\{2,3\},\ldots,\{k,1\})$ be the unique cycle of $G$, where $k \geq 4$ is even. Let $\{v_{1},\ldots,v_{k},w_{1},\ldots,w_{k}\}$ be the set of vertices of $\overline{C}$, where $\{v_{1},\ldots,v_{k}\} \cap \{w_{1},\ldots,w_{k}\}=\emptyset$. Every odd cycle in $H$ contains at least one of the edges $z_{1k}=\{v_{1},w_{k}\}$ or $z_{k1}=\{w_{1},v_{k}\}$, since the subgraph $F$ of $H$ with the edge set $E(F)=E(\overline{C}) \setminus \{z_{1k},z_{k1}\}$ is bipartite and a partition of its vertices is $$\{v_{1},v_{3},\ldots,v_{k-1},w_{2},w_{4},\ldots,w_{k}\} \cup \{v_{2},v_{4},\ldots,v_{k},w_{1},w_{3},\ldots,w_{k-1}\}.$$ We will prove that any two odd cycles in $H$ have at least 2 vertices in common. Suppose that there exist two odd cycles $\Gamma_{1}$ and $\Gamma_{2}$ in $H$ which share at most one vertex. Let, say, that $\Gamma_{1}$ contains $z_{1k}$ and $\Gamma_{2}$ contains $z_{k1}$. We can take the cycles $\Gamma_{1}$ and $\Gamma_{2}$ to start from the vertices $v_{1}$ and $w_{1}$, respectively. Moreover, we can assume that $z_{1k}$ and $z_{k1}$ are the first edges of $\Gamma_{1}$ and $\Gamma_{2}$, respectively. We claim that the second edge of $\Gamma_{1}$ is $\{w_{k},w_{k-1}\}$. If $\{w_{k},w_{k-1}\}$ is not the second edge, then either $\{w_{k},v_{k}\}$ is the second edge of $\Gamma_{1}$ or there exists a vertex $q \in T_{i}^{\star}$ with $q \neq v_{k}$ such that $\{w_{k},q\}$ is the second edge of $\Gamma_{1}$. In the latter case there exists a path in $\Gamma_{1}$ of length $>2$ connecting $w_{k}$ with $v_{k}$. Since $\Gamma_{1}$ is a cycle, we have that, in both cases, $\{v_{k},v_{k-1}\}$ is an edge of $\Gamma_{1}$. Also $\{w_{1},v_{k}\}$ is the first edge of $\Gamma_{2}$, so either $\{v_{k},v_{k-1}\}$ is the second edge of $\Gamma_{2}$ or there exists a path in $\Gamma_{2}$ of length $\geq$ 1 connecting $v_{k}$ with $w_{k}$. In both cases we arrive at a contradiction, since the cycles $\Gamma_{1}$, $\Gamma_{2}$ have at most one common vertex. Consequently, $\{w_{k},w_{k-1}\}$ is the second edge of $\Gamma_{1}$ and analogously we have that $\{v_{k},v_{k-1}\}$ is the second edge of $\Gamma_{2}$. Using similar arguments we conclude that $\{v_{1},w_{k}\}, \{w_{k},w_{k-1}\}, \ldots,\{w_{3},w_{2}\}$ are all edges of $\Gamma_{1}$, while $\{w_{1},v_{k}\}, \{v_{k},v_{k-1}\}, \ldots,\{v_{3},v_{2}\}$ are all edges of $\Gamma_{2}$. We claim that $\{w_{2},w_{1}\}$ is also an edge of $\Gamma_{1}$. Suppose not, then either $\{w_{2},v_{2}\}$ is an edge of $\Gamma_{1}$ or there exists a vertex $q' \in T_{j}^{\star}$ with $q' \neq w_{1}$ such that $\{w_{2},q'\}$ is an edge of $\Gamma_{1}$. In both cases there exists a path in $\Gamma_{1}$ of length $\geq 1$ connecting $w_{2}$ with $v_{2}$. Thus $\{v_{2},v_{1}\}$ is an edge of $\Gamma_{1}$. But $\{v_{3},v_{2}\}$ is an edge of $\Gamma_{2}$, so either $\{v_{2},v_{1}\}$ is an edge of $\Gamma_{2}$ or there exists a path in $\Gamma_{2}$ of length $\geq 1$ connecting $v_{2}$ with $w_{2}$. Since the cycles $\Gamma_{1}$, $\Gamma_{2}$ have at most one vertex in common, we arrive at a contradiction. Thus $\{w_{2},w_{1}\}$ is an edge of $\Gamma_{1}$ and analogously $\{v_{2},v_{1}\}$ is an edge of $\Gamma_{2}$. But then $\Gamma_{1}$, $\Gamma_{2}$ have 2 vertices in common, namely $v_{1}$ and $w_{1}$, contradicting our assumption. As a consequence any two odd cycles in $H$ have at least 2 vertices in common. From Theorem \ref{Circuits of a graph} we have that every circuit of $I_{H}$ is of the form $B_{\Gamma}$, where $\Gamma$ is an even cycle of $H$. Consequently, the universal Gr\"{o}bner basis of $P_{G}$ consists of all binomials of the form $B_{\Gamma}$, where $\Gamma$ is an even cycle of $H$. \hfill $\square$

\begin{rem1} \label{Upper bound sharp} {\rm Let $C=(\{1,2\},\{2,3\},\ldots,\{k,1\})$ be an even cycle of $G$ of length $k$, then $P_{C}=I_{\overline{C}}$, and therefore, from Remark \ref{Odd Cycles Basic}, the maximum degree of a binomial in the universal Gr\"{o}bner basis of $P_{C}$ is equal to $k$. Notice that $\left \lfloor{ \frac{k+k+1}{2}} \right \rfloor=k$.}
\end{rem1}

\begin{prop1} Let $G$ be a simple graph on the vertex set $\{1,\ldots,n\}$ and $C=(\{p_{1},p_{2}\}, \{p_{2},p_{3}\},\ldots,\{p_{k},p_{1}\})$ be an odd cycle of $G$ of length $k \geq 3$. Consider the graph $C^{\star}$ on the vertex set $\{p_{1},\ldots,p_{k}, q_{1},\ldots,q_{k}\}$, where $\{p_{1},\ldots,p_{k}\} \cap \{q_{1},\ldots,q_{k}\}=\emptyset$, and let $C'$ be the odd cycle $$C'=(\{q_{1},q_{2}\},\{q_{2},q_{3}\},\ldots,\{q_{k-1},q_{k}\},\{q_{k},q_{1}\}).$$ Then a binomial $B \in P_{C}$ belongs to the universal Gr\"{o}bner basis of $P_{C}$ if and only if $B=B_{\Gamma}$ where $\Gamma$ is one of the following even closed walks of $C^{\star}$: \begin{enumerate} \item $\Gamma$ is an even cycle of $C^{\star}$, \item $\Gamma=(C, \{p_{i},q_{i}\}, C')$ where $1 \leq i \leq k$.
\end{enumerate}
\end{prop1}

\noindent \textbf{Proof.} Since the equality $P_{C}=I_{C^{\star}}$ holds, we deduce from Proposition \ref{Proposit1} that the toric ideal $I_{C^{\star}}$ is complete intersection. Thus, from Proposition 6.1 in \cite{BGR}, the graph $C^{\star}$ has exactly two vertex disjoint odd cycles, namely $C$ and $C'$. By Lemma 3.2 in \cite{OH3}, every primitive binomial $B \in I_{C^{\star}}$ is of the form $B=B_{\Gamma}$, where $\Gamma$ is one of the following even closed walks:\\ (1) $\Gamma$ is an even cycle of $C^{\star}$,\\ (2) $\Gamma$ consists of two odd cycles of $C^{\star}$ intersecting in exactly one vertex,\\ (3) $\Gamma=(C, \{p_{i},q_{i}\}, C')$ where $1 \leq i \leq k$, i.e. it consists of the vertex disjoint odd cycles $C$, $C'$ joined by the edge $\{p_{i},q_{i}\}$.\\ But $C^{\star}$ has no vertex of degree greater than three, so there are no two odd cycles of $C^{\star}$ intersecting in exactly one vertex. Furthermore every primitive binomial of $I_{C^{\star}}$ is also a circuit. Consequently the universal Gr\"{o}bner basis of $P_{C}$ consists of all binomials of the form $B_{\Gamma}$, where either $\Gamma$ is an even cycle of $C^{\star}$ or $\Gamma=(C, \{p_{i},q_{i}\}, C')$. \hfill $\square$

\begin{rem1} {\rm Let $G$ be a connected non-bipartite graph with exactly one cycle. As Example \ref{Graver} demonstrates, the equality $U(P_{G})=\mathcal{C}(I_{G^{\star}})$ does not hold in general. From Theorem 3.4 in \cite{Tata-Tho} we have that a binomial $B$ belongs to the universal Gr\"{o}bner basis of $P_{G}$ if and only if $B=B_{\Gamma}$, where $\Gamma$ is a mixed even closed walk of $G^{\star}$. For the definition of a mixed walk see \cite{Tata-Tho}.}
\end{rem1}

Another interesting case of a bipartite graph is that of a tree $G$. The next theorem provides a graph theoretic characterization of the elements in the universal Gr\"{o}bner basis of $P_{G}$.

\begin{thm1} \label{UniversGrobn} Let $G$ be a tree on the vertex set $\{1,\ldots,n\}$, then the universal Gr\"{o}bner basis of $P_{G}$ is given by the set $$U(P_{G})=\{B_{\Gamma}|\Gamma \ \textrm{is an even cycle of} \ G^{\star}\}.$$
\end{thm1}

\noindent \textbf{Proof.} By the proof of Lemma \ref{Tree}, the graph $G^{\star}$ is bipartite and also $P_{G}=I_{G^{\star}}$. So $I_{A_{G}}=I_{G^{\star}}$, and therefore we have, from Theorem \ref{UGB}, that $U(P_{G})=\mathcal{C}(I_{G^{\star}})$. Using Theorem \ref{Circuits of a graph} we conclude that the universal Gr\"{o}bner basis of $P_{G}$ consists of all the binomials $B_{\Gamma}$, where $\Gamma$ is an even cycle of $G^{\star}$. \hfill $\square$

\begin{rem1} {\rm Given a tree $G$, we have that the cardinality of $U(P_{G})$ is equal to the number of paths in $G$. Furthermore the maximum degree of a binomial in the universal Gr\"{o}bner basis of $P_{G}$ equals $l+1$, where $l$ is the length of the longest path of $G$. In particular, for a star graph $G$ on the vertex set $\{1,\ldots,n\}$ with $n \geq 3$, the universal Gr\"{o}bner basis of $P_{G}$ consists of $\frac{n(n-1)}{2}$ binomials and the maximum degree of a binomial in $U(P_{G})$ is equal to three.}
\end{rem1}

\begin{rem1} \label{Upper bound is sharp} {\rm For a path graph $G$ with $n$ vertices, we have that the maximum degree of a binomial in the universal Gr\"{o}bner basis of $P_{G}$ equals $n=\left \lfloor{ \frac{n+(n-1)+1}{2}} \right \rfloor$.}
\end{rem1}

\end{document}